\documentclass[11pt]{article}
\usepackage{latexsym,amsmath,color,amsthm,amssymb,epsfig,graphicx,mathrsfs}
\usepackage{graphicx}
\usepackage{amssymb}
\usepackage[left=1in,top=1in,right=1in,bottom=1in]{geometry}
\usepackage[linktocpage=true]{hyperref}
\usepackage{setspace}
\usepackage{amssymb, amsmath, amsthm, graphicx,mathrsfs}

\def\qed{\hfill\ifhmode\unskip\nobreak\fi\quad\ifmmode\Box\else\hfill$\Box$\fi}
\def\ite#1{\hfill\break${}$\hbox to 50pt {\quad(#1)\hfill}}
\newtheorem{thm}{Theorem}[section]
\newtheorem{cor}[thm]{Corollary}

\newtheorem{lem}[thm]{Lemma}

\begin{document}

\title{\vspace{-0.5in}The maximum number of cliques in graphs without long cycles}

\author{
{{Ruth Luo}}\thanks{University of Illinois at Urbana--Champaign, Urbana, IL 61801, USA. E-mail: {\tt ruthluo2@illinois.edu}. Research of this author
is supported in part by NSF grant DMS-1600592.}}

\maketitle

\vspace{-0.3in}

\begin{abstract}
The Erd\H{o}s--Gallai Theorem states that for $k\geq 3$ every graph on $n$ vertices with more than $\frac{1}{2}(k-1)(n-1)$ edges contains a cycle of length at least $k$. Kopylov proved a strengthening of this result for 2-connected graphs with extremal examples $H_{n,k,t}$ and $H_{n,k,2}$. In this note, we generalize the result of Kopylov to bound the number of $s$-cliques in a graph with circumference less than $k$. Furthermore, we show that the same extremal examples that maximize the number of edges also maximize the number of cliques of any fixed size. Finally, we obtain the extremal number of $s$-cliques in a graph with no path on $k$-vertices.

\end{abstract}

\section{Introduction} In \cite{ErdGal59}, Erd\H{o}s and Gallai determined $ex(n, P_k)$, the maximum number of edges in an $n$-vertex graph that does not contain a copy of the path on $k$ vertices, $P_k$. This result was a corollary of the following theorem:

\begin{thm}[Erd\H{o}s and Gallai~\cite{ErdGal59}]\label{ErdGallaiCyc}
Let $G$ be an $n$-vertex graph with more than $\frac{1}{2}(k-1)(n-1)$ edges, $k \geq 3$.
Then $G$ contains a cycle of length at least $k$.
\end{thm}

To obtain the result for paths, suppose $G$ is an $n$-vertex graph with no copy of $P_k$. Add a new vertex $v$ adjacent to all vertices in $G$, and let this new graph be $G'$. Then $G'$ is an $n+1$-vertex graph with no cycle of length $k+1$ or longer, and so $e(G)+ n = e(G') \leq \frac{1}{2}kn$ edges.

\begin{cor}[Erd\H{o}s and Gallai~\cite{ErdGal59}]\label{ErdGallaiPath}
Let $G$ be an $n$-vertex graph with more than $\frac{1}{2}(k-2)n$ edges, $k\geq 2$.
Then $G$ contains a copy of $P_k$.
\end{cor}

Both results are sharp with the following extremal examples: for Theorem \ref{ErdGallaiCyc}, when $k-2$ divides $n-1$, take any connected $n$-vertex graph whose blocks (maximal connected subgraphs with no cut vertices) are cliques of order $k-1$. For Corollary \ref{ErdGallaiPath}, when $k -1$ divides $n-1$, take the $n$-vertex graph whose connected components are cliques of order $k - 1$.

There have been several alternate proofs and sharpenings of the Erd\H{o}s-Gallai theorem
including results by Woodall \cite{Woodall}, Lewin \cite{Lewin}, Faudree and Schelp\cite{FaudScheB,FaudSche75}, and Kopylov \cite{Kopy} -- see~\cite{FS224} for further details.

The strongest version was that of Kopylov who improved the Erd\H{o}s--Gallai bound for 2-connected graphs. To state the theorem, we first introduce the family of extremal graphs. 

Fix $k\geq 4$, $n \geq k$, $\frac{k}{2} > a\geq 1$. Define the $n$-vertex graph $H_{n,k,a}$ as follows.
 The vertex set of $H_{n,k,a}$ is partitioned into three sets $A,B,C$ such that $|A| = a$, $|B| = n - k + a$ and $|C| = k - 2a$
 and the edge set of $H_{n,k,a}$ consists of all edges between $A$ and $B$ together with all edges in $A \cup C$.

Note that when $a \geq 2$, $H_{n,k,a}$ is 2-connected, has no cycle of length $k$ or longer, and $e(H_{n,k,a}) = {k-a \choose 2} + (n-k+a)a$.

\begin{figure}[!ht]
  \centering
    \includegraphics[width=0.25\textwidth]{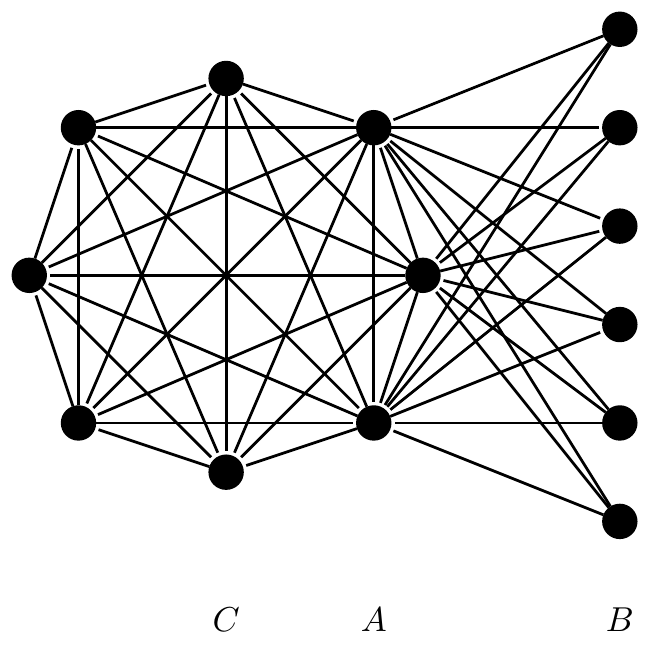}
  \caption{$H_{14,11,3}$}
\end{figure}

{\bf Definition}. Let $f_s(n,k, a):= {k - a \choose s} + (n-k + a){a \choose s-1}$, where $f_2(n,k,a) = e(H_{n,k,a})$.

By considering the second derivative, one can check that $f_s(n,k, a)$ is convex in $a$ in the domain $[1, \lfloor (k-1)/2 \rfloor]$, thus it attains its maximum at one of the endpoints $a = 1$ or $a =   \lfloor (k-1)/2 \rfloor$. 
\begin{thm}[Kopylov~\cite{Kopy}]\label{kopycycle} Let $n \geq k \geq 5$ and let $t = \lfloor \frac{k-1}{2}\rfloor$. If $G$ is a 2-connected $n$-vertex graph with \[e(G) \geq \max\{f_2(n, k, 2), f_2(n, k, t)\},\]
then either $G$ has a cycle of length at least $k$, or $G = H_{n,k, 2}$, or $G = H_{n,k, t}$. 
\end{thm}

 \begin{figure}[!ht]
  \centering
    \includegraphics[height=0.15\textwidth]{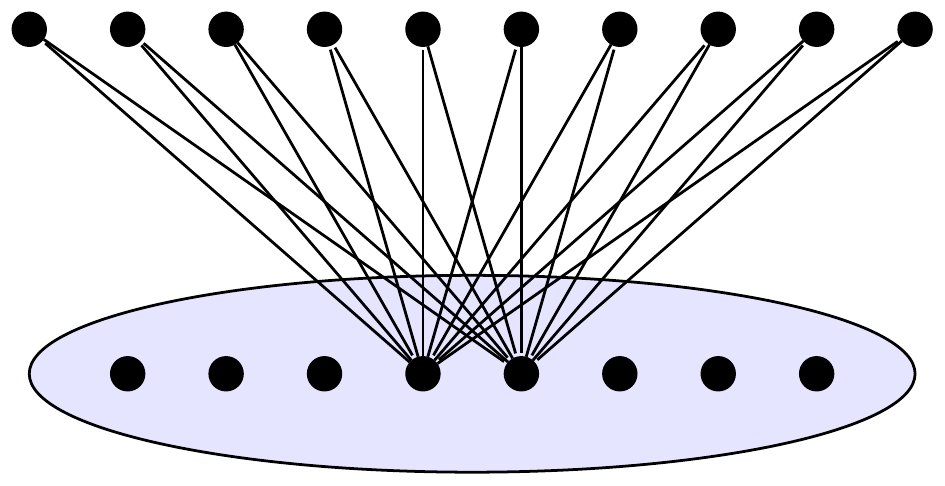}\;\;\;
    \includegraphics[height=0.15\textwidth]{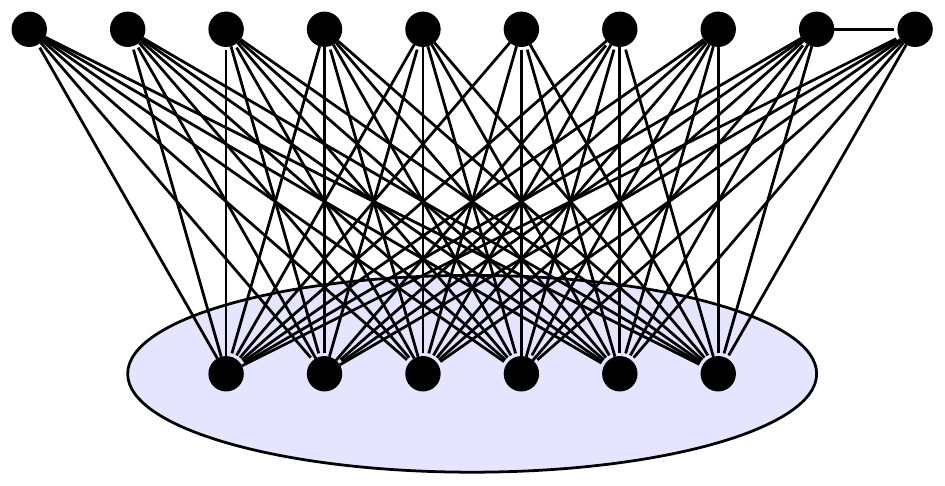}\;\;\;
    \includegraphics[height=0.15\textwidth]{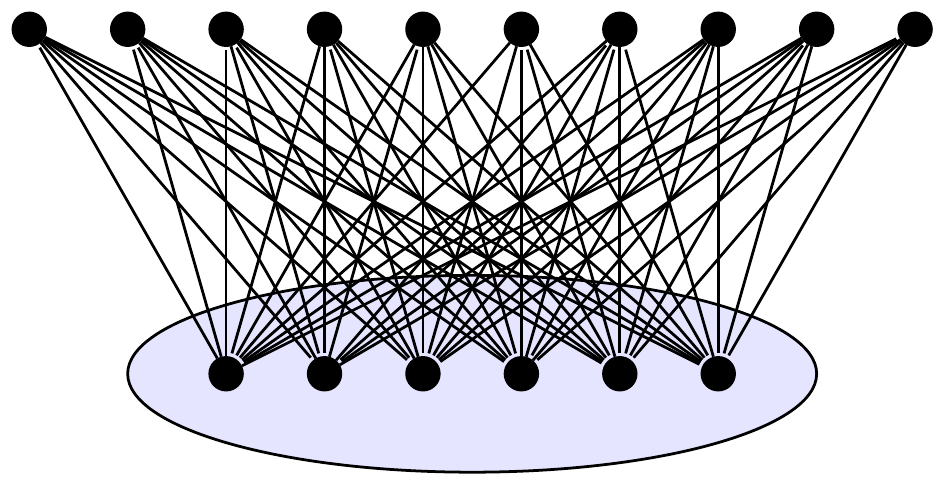}
  \caption{$ H_{n,k,2}, H_{n,k,t} (k = 2t+1), H_{n,k,t} (k = 2t + 2)$; ovals denote complete subgraphs of sizes $k-2$, $t$, and $t$, respectively.}
\end{figure}

It is straight-forward to check that any 2-connected graph that is not a triangle has a cycle of length $4$ or greater, and so the theorem covers all nontrivial choices of $k$.  This theorem also implies Theorem \ref{ErdGallaiCyc} by applying induction to each block of the graph. 

We consider a generalized Tur\'{a}n-type problem. Fix graphs $T$ and $H$, and define the function $ex(n,T,H)$ to be the maximum number of (unlabeled) copies of $T$ in an $H$-free graph on $n$ vertices. When $T = K_2$, we have the usual extremal number $ex(n,T,H) = ex(n, H)$.

There are many notable papers studying the $ex(n,T,H)$ function for different combinations of $T$ and $H$. Erd\H{o}s~\cite{Erd} proved that for $s \leq r$, among all $n$-vertex graphs that forbid $K_{r+1}$, the Tur\'an graph  (i.e., the balanced complete $r$-partite graph) maximizes the number of copies of $K_s$.  Hatami, Hladk\'y, Kr\'al', Norine, and Razborov~\cite{raz} and independently Grzesik~\cite{grz} proved $ex(n, C_5, K_3) = (n/5)^5$ whenever $n$ is divisible by 5 using the method of flag algebras. On the other hand, Bollob\'as and Gy\H{o}ri~\cite{bol} proved $(1+o(1)) \frac{1}{3\sqrt{3}}n^{3/2} \leq ex(n, K_3, C_5) \leq (1+o(1)) \frac{5}{4}n^{3/2}$, and later Gy\H{o}ri and Li~\cite{GL} proved an upper bound for $ex(n, K_3, C_{2k+1})$ in terms of $ex(n, C_{2k})$. This bound was improved by F\"uredi and \"Ozkahya~\cite{furedi} and then later improved again by Alon and Shikhelman~\cite{alonsh}. In the same paper, Alon and Shikhelman proved $ex(n, K_s, K_{r,t}) = \Theta(n^{s-{s \choose 2}/r})$ for certain values of $r$, $s$, and $t$, among other results. 

Furthermore, such generalized Tur\'an-type results for graphs can be instrumental for proving related extremal results in hypergraphs. For example, F\"uredi and \"Ozkahya~\cite{furedi} used their upper bounds for the number of triangles in graphs without cycles of fixed lengths to give an upper bound for the number of hyperedges in 3-uniform hypergraphs without Berge-cycles of a fixed length. 



In this note, we give an upper bound for the number of $s$-cliques in a graph without cycles of length $k$ or greater (i.e., circumference less than $k$). We also obtain $ex(n,K_s, P_k)$.

{\bf Definition}. For $s \geq 2$, let $N_s(G)$ denote the number of unlabeled copies of $K_s$ in $G$, e.g., $N_2(G) = e(G)$.

Our main result is a generalization of Kopylov's result, Theorem \ref{kopycycle}. In particular, we show that the same extremal examples that maximize the number of edges among $n$-vertex 2-connected graphs with circumference less than $k$ also maximize the number of cliques of any size. Our main results are the following:

\begin{thm}\label{main} Let $n \geq k \geq 5$ and let $t = \lfloor \frac{k-1}{2}\rfloor$. If $G$ is a 2-connected $n$-vertex graph with circumference less than $k$, then \[N_s(G) \leq \max\{f_s(n, k, 2), f_s(n, k, t)\}.\]
\end{thm}
Again, this theorem is sharp with the same extremal examples $H_{n,k, 2}$ and $H_{n,k,t}$.

This theorem implies the cliques version of Theorem \ref{ErdGallaiCyc}:

\begin{cor}\label{newcycle}Let $n \geq k \geq 4$. If $G$ is an $n$-vertex graph with circumference less than $k$, then \[N_s(G) \leq \frac{n-1}{k-2} {k-1 \choose s}.\]
\end{cor}

Unlike the edges case, Theorem \ref{main} unfortunately does not easily imply $ex(n, K_s, P_k)$. However, a Kopylov-style argument very similar to the proof of Theorem \ref{main} gives the result for paths.

\begin{thm} \label{newpath1}
Let $n \geq k\geq 4$ and let $G$ be an $n$-vertex connected graph with no path on $k$ vertices. Let $t = \lfloor (k - 2)/2 \rfloor$. Then $N_s(G) \leq \max\{f_s(n, k-1, 1), f_s(n, k-1, t)\}$. 
\end{thm}
We have sharpness examples $H_{n, k-1, 1}$ and $H_{n,k-1, t}$. Finally, using induction on the number of components gives the following result:

\begin{cor} \label{newpath2}
$ex(n, K_s, P_k) = \frac{n}{k - 1} {k - 1 \choose s}$. 
\end{cor}
And the same extremal examples as for Corollary \ref{ErdGallaiPath} apply.

The proofs for Corollary \ref{newcycle}, Theorem \ref{newpath1}, and Theorem \ref{newpath2} are given in Section 3 of this paper. We first prove Theorem \ref{main}.

\section{Proof of Theorem \ref{main}}

Let $G$ be an edge-maximal counterexample. Then $G$ is {\em $k$-closed}, i.e., adding any additional edge to $G$ creates a cycle of length at least $k$. In particular, for any nonadjacent vertices $x$ and $y$ of $G$, there exists a path of at least $k-1$ edges between $x$ and $y$. We will use the following lemma:

\begin{lem}[Kopylov \cite{Kopy}]\label{path} Let $G$ be a 2-connected $n$-vertex graph with a path $P$ of $m$ edges with endpoints $x$ and $y$. For $v \in V(G)$, let $d_P(v) = |N(v) \cap V(P)|$. Then $G$ contains a cycle of length at least $\min\{m+1, d_P(x) + d_P(y)\}.$
\end{lem}

Our first goal is to show that $G$ contains a large ``core'', i.e.,  a subgraph with large minimum degree. 
For this, we use the notion of \emph{disintegration}.

{\bf Definition}: For a natural number $\alpha$ and a graph $G$, the \emph{$\alpha$-disintegration} of 
a graph $G$ is the  process of iteratively removing from $G$ the vertices with degree 
at most $\alpha$  until the resulting graph has minimum degree at least $\alpha + 1$ or is empty. 
This resulting subgraph $H = H(G, \alpha)$ will be called the $(\alpha+1)$-{\em core} of $G$. It is well known
that $H(G, \alpha)$ is unique and does not depend on the order of vertex deletion (for instance, see \cite{core}).

Let $H(G, t)$ denote the $(t+1)$-core of $G$, i.e., the resulting graph of applying $t$-disintegration to $G$. We claim that
\begin{equation*}
\mbox{\em
 $H(G,t)$ is nonempty.}
\end{equation*}

Suppose $H(G,t)$ is empty. In the disintegration process, every time a vertex of degree at most $t$ is removed, we delete at most ${t \choose s-1}$ copies of $K_s$. For the last $\ell \leq t$ vertices, we remove at most ${\ell-1 \choose s-1}$ copies of $K_s$ with each deletion. Thus
\begin{eqnarray*}
N_s(G) &\leq& (n-t){t \choose s-1} + {t-1 \choose s-1} + {t-2 \choose s-1} + \ldots + {0 \choose s-1}\\
& =& (n-t){t \choose s-1} + {t \choose s} \\
& = & (n-(t+1)){t \choose s-1} + {t+1 \choose s}\\
& \leq &  f_s(n,k, t),
\end{eqnarray*}

a contradiction.

Therefore $H(G,t)$ is nonempty. Next we show that 
\begin{equation*}
\mbox{\em
 $H(G,t)$ is a complete graph.}
\end{equation*}
If there exists a nonedge of $H(G,t)$, then in $G$, there is a path of length at least $k-1$ edges with these vertices as its endpoints. Among all nonadjacent pairs of vertices in $H(G,t)$, choose $x,y$ such that there is a longest path $P$ in $G$ with endpoints $x$ and $y$. By maximality of $P$, all neighbors of $x$ in $H(G,t)$ lie in $P$: if $x$ has a neighbor $x' \in H(G,t) - P$, then either $x'y \in E(G)$ and $x'P$ is a cycle of length at least $k$, or $x'y \notin E(G)$ and so $x'P$ is a longer path. Similar for $y$. Hence, by Lemma~\ref{path}, $G$ has a cycle of length at least $\min\{k, d_P(x) + d_P(y)\} = \min\{k, 2(t+1)\} = k$, a contradiction.  

Now let $r = |V(H(G,t)|$. Each vertex in $H(G,t)$ has degree at least $t+1$, so $r \geq t+2$. Also, if $r \geq k -1$, as $G$ is 2-connected and $H(G,t)$ is a clique, we can extend a path on $r$ vertices of $H(G,t)$ to a cycle of length at least $r+1 \geq k$, a contradiction. Therefore $t+2 \leq r \leq k - 2$. In particular, $2 \leq k - r \leq t$. Apply $(k - r)$-disintegration to $G$, and let $H(G, k - r)$ be the resulting graph. Then $H(G, t) \subseteq H(G, k - r)$. 

If $H(G, t) = H(G, k - r)$, then \[N_s(G)  \leq {r \choose s} + (n-r){k - r \choose s-1} = f_s(n, k, k - r) \leq \max\{f_s(n,k, 2),f_s(n,k,t)\}\] by the convexity of $f_s$. 
Therefore, $H(G, t)$ is a proper subgraph of $H(G, k - r)$, and there must be a nonedge between a vertex in $H(G,t)$ and a vertex in $H(G, k - r)$. Among all such pairs, choose $x \in H(G,t)$ and $y \in H(G, k - r)$ to have a longest path $P$ between them. As before, $P$ contains at least $k-1$ edges, and each neighbor of $x$ in $H(G,t)$ and each neighbor of $y$ in $H(G, k-r)$ lie in $P$. Then $G$ contains a cycle of length at least $\min\{k,  (r-1) + (k - r + 1)\} = k$, a contradiction. \qed

\section{Proof of Corollary \ref{newcycle}, Theorem \ref{newpath1}, and Corollary \ref{newpath2}}

{\em Proof of Corollary \ref{newcycle}}.
Define $g_s(n, k) = \frac{n-1}{k-2}{k - 1 \choose s}$ and $t = \lfloor \frac{k - 1}{2}\rfloor$. One can check that when $n \geq k$, \[g_s(n, k) \geq \max\{f_s(n,k, t), f_s(n, k, 2)\}.\] 

 Fix a graph $G$ on $n$ vertices with circumference less than $k$. If $G$ is disconnected, simply apply induction to each component of $G$ to obtain the desired result. Therefore we may assume $G$ is connected. We induct on the number of blocks of $G$. First suppose $k \geq 5$. If $G$ is a block, i.e., 2-connected,  then either $n \leq k - 1$, and so $N_s(G) \leq {|V(G)| \choose s} \leq g_s(n, k)$, or $n \geq k$, and so by Theorem \ref{main}, $N_s(G) \leq \max\{f_s(n,k, t), f_s(n, k, 2)\} \leq g_s(n,k)$. 
 
Otherwise, consider the {\em block-cut tree} of $G$---the tree whose vertices correspond to blocks of $G$ such that two vertices in the tree are adjacent if and only if the corresponding blocks in $G$ share a vertex. Let $B_1$ be a block in $G$ corresponding to a leaf-vertex in the block-cut tree such that $B_1$ and its complement are connected by the cut vertex $v$. Set $B_2 = G - B_1 + \{v\}$.  Apply the induction hypothesis to $B_1$ and $B_2$ to obtain
\begin{eqnarray*}
N_s(G) &=& N_s(B_1) + N_s(B_2) \leq g_s(|B_1|, k) + g_s(n-|B_1|+1, k)\\
&=& \frac{|B_1|-1}{k-2}{k-1 \choose s} + \frac{(n-|B_1|+1) - 1}{k-2}{k-1 \choose s}\\
&=& g_s(n, k).\end{eqnarray*} 

If $k=4$, then either $G$ is a forest or $G$ has circumference 3. In the second case, each block of $G$ is either a triangle or an edge. Thus $N_s(G) \leq g_s(n,k)$ in both cases. 
\qed 
\\ \\
The proof of Theorem \ref{newpath1} follows the same steps as the proof of Theorem \ref{main}. As some details here will be omitted to prevent repetition, it is advised that the reader first reads the proof of Theorem \ref{main}.

{\em Proof of Theorem \ref{newpath1}}. Suppose for contradiction that $N_s(G) > \max\{f_s(n, k-1, 1),f_s(n, k-1, t)\}$ where $t = \lfloor (k-2)/2\rfloor$.  Let $G_0$ be the graph obtained by adding a dominating vertex $v_0$ adjacent to all of $V(G)$. Then $G_0$ is 2-connected, has $n+1$ vertices, and contains no cycle of length $k + 1$ or greater. Let $G'$ be the $k +1$-closure of $G_0$ (i.e., add edges to $G_0$ until any additional edge creates a cycle of length at least $k +1$). Denote by $N'_s(G')$ the number of $K_s$'s in $G'$ that do not contain $v_0$. Thus $N'_s(G') \geq N'_s(G_0) = N_s(G)$.  Apply $(t+1)$-disintegration to $G'$, where if necessary, we delete $v_0$ last. Let $H(G', t+1)$ be the resulting graph of the disintegration. If $H(G', t+1)$ is empty, then at the time of deletion each vertex has at most $t$ neighbors that are not $v_0$. Hence
\[N'_s(G') \leq  (n-(t+1)){t \choose s-1} + {t+1 \choose s} \leq f_s(n, k - 1, t),\]
a contradiction. 

The same argument as in the proof of Theorem \ref{main} also shows that $H(G', t+1)$ is a complete graph, otherwise there would be a cycle of length at least $2(t+2) \geq (k -1)+2$ in $G'$. Note that $v_0$ must be contained in $H(G', t+1)$ as it is adjacent to all vertices in $G'$. Set $|V(H(G', t+1))| = r$ where $t+3 \leq r \leq k - 1$ (and so $k - r \geq 1$). In particular, $(k + 1) - r \leq t+1$. Apply $(k + 1 - r)$-disintegration to $G'$. If $H(G', t+1) \neq H(G', k + 1 - r)$, then again we can find a cycle of length at least $(r-1) + k + 2 - r = k +1$. Otherwise, suppose $H(G', t+1) = H(G', k + 1 - r)$. In $H(G', t+1)$, the number of $s$-cliques that do not include $v_0$ is ${r-1 \choose s}$, and in $V(G) - V(H(G', k + 1 - r))$, every vertex had at most $k - r$ neighbors that were not $v_0$ at the time of its deletion. We have
\begin{eqnarray*}
N'_s(G') &\leq& {r-1 \choose s} + (n+1 - r){ k - r \choose s-1}\\
& =& f_s(n, k - 1, k - r) 
 \leq \max\{f_s(n, k - 1, 1),f_s(n,k-1,t)\},\end{eqnarray*}
a contradiction. 
\qed

{\em Proof of Corollary \ref{newpath2}}.
Define $h_s(n, k) = \frac{n}{k -1} {k - 1 \choose s}$, and note that when $n \geq k$, \[h_s(n, k) \geq \max\{f_s(n, k-1, t), f_s(n, k-1, 1)\}.\] 

We induct on the number of components in $G$. First suppose $k \geq 4$. If $G$ is connected, then either $n \leq k - 1$, in which case $N_s(G) \leq {|V(G)| \choose s} \leq h_s(n,k)$, or $n \geq k$ and $N_s(G) \leq \max\{ f_s(n, k-1, 1),f_s(n, k-1, t)\} \leq h_s(n, k)$.   Otherwise if $G$ is not connected, let $C_1$ be a component of $G$. Then $N_s(G) = N_s(C_1) + N_s(G - C_1) \leq h_s(|C_1|, k) + h_s(n-|C_1|, k) = h_s(n,k)$. 

If $k=3$ (the cases $k \leq 2$ are not interesting), then the longest path in $G$ has two vertices. It follows that $G$ is the union of a matching and isolated vertices. Therefore $N_s(G) \leq h_s(n,k)$.
\qed

\vspace{12mm}
{\bf Acknowledgment.} The author would like to thank Alexandr Kostochka, Zolt\'an F\"uredi, and Jacques Verstra\"ete for their guidance and for sharing their knowledge on this topic.


\begin{thebibliography}{99}




\bibitem{alonsh} %
N. Alon and C. Shikhelman,
Many $T$ copies in $H$-free graphs,
J. Combin. Theory Ser. B. {\bf 121} (2016), 146--172.

\bibitem{bol} B. Bollob\'as and E. Gy\H{o}ri, Pentagons vs. triangles, Discrete Math. {\bf 308} (2008), 4332–4336

\bibitem{Erd}P. Erd\H{o}s, On the number of complete subgraphs contained in certain graphs, Magyar Tud. Akad.
Mat. Kut. Int. K˝ozl, {\bf 7} (1962), 459-–474.

\bibitem{ErdGal59}
P. Erd\H{o}s and T. Gallai,
On maximal paths and circuits of graphs,
Acta Math. Acad. Sci. Hungar. {\bf 10} (1959), 337--356.

\bibitem{FaudScheB}
R. J. Faudree and R. H. Schelp, Ramsey type results,
Infinite and Finite Sets, {Colloq. Math. J. Bolyai} {\bf 10}, (ed. A. Hajnal et al.), North-Holland,
Amsterdam, 1975, pp. 657--665.

\bibitem{FaudSche75}
R. J. Faudree and R. H. Schelp,
Path Ramsey numbers in multicolorings,
J. Combin. Theory Ser. B. {\bf 19} (1975), 150--160.

\bibitem{furedi} Z. F\"uredi and L. \"Ozkahya, On 3-uniform hypergraphs without a cycle of a given length, Discrete Applied Mathematics
{\bf 216}, Part 3, (2017), 582-–588.


\bibitem{FS224}
Z. F\"uredi and M. Simonovits,
  The history of degenerate (bipartite) extremal graph problems,
Bolyai Math. Studies {\bf 25} pp. 169--264,
Erd\H{o}s Centennial (L. Lov\'asz, I. Ruzsa, and V. T. S\'os, Eds.) Springer, 2013.
Also see: {\tt arXiv:1306.5167}.


\bibitem{grz}
A. Grzesik, On the maximum number of five-cycles in a triangle-free graph, J. Combin. Theory Ser. B. {\bf 102.5} (2012),  1061-1066.
\bibitem{GL} E. Gy\H{o}ri and H. Li, The maximum number of triangles in $C_{2k+1}$-free graphs, Combinatorics,
Probability and Computing {\bf21}(1-2), (2012), 187-191.

\bibitem{raz} %
H. Hatami, J. Hladk\'y, D. Kr\'al', S. Norine, and A. Razborov, On the number of pentagons in triangle-free graphs, J. Combin. Theory Ser. A. {\bf 120} (2013) no. 3, 722--732.

\bibitem{Kopy}
G. N. Kopylov,
Maximal paths and cycles in a graph,
Dokl. Akad. Nauk SSSR {\bf 234} (1977),  19--21.
(English translation: Soviet Math. Dokl. {\bf 18} (1977), no. 3, 593--596.)


\bibitem{Lewin}
M. Lewin,
On maximal circuits in directed graphs,
J. Combin. Theory Ser. B. {\bf 18} (1975), 175--179.


\bibitem{core} B. Pittel, J. Spencer, and N. Wormald, Sudden emergence of a giant $k$-core in a random graph, J. Combin. Theory Ser. B. {\bf 67} (1996), 111–151.

\bibitem{Woodall}
D. R. Woodall,
Maximal circuits of graphs I,
Acta Math. Acad. Sci. Hungar. {\bf 28} (1976), 77--80.


\end{thebibliography}
\end{document}